\documentclass[12pt,a4paper,twoside]{article}
\usepackage[T2A]{fontenc}
\usepackage[english,russian]{babel}

\fontsize{12}{21} \selectfont

\newtheorem{theorem}{Theorem}[section]

\newtheorem{definition}{Definition}[section]

\begin{document}

\noindent

\makeatletter
\renewcommand{\@evenhead}{\hfil
{\bf Andrei I. Bodrenko.}}
\renewcommand{\@oddhead}{\hfil
\small{\bf \underline{The solution of the Minkowski problem for
closed surfaces in Riemannian space.}} }

\noindent
\bigskip
\begin{center}
{\large \bf The solution of the Minkowski problem for closed
surfaces in Riemannian space.}
\end{center}
\medskip

\begin{center}
{\bf Andrei~I.~Bodrenko} \footnote{\copyright  Andrei~I.~Bodrenko,
associate professor,
Department of Mathematics,\\
Volgograd State University,
University Prospekt 100, Volgograd, 400062, RUSSIA.\\
E.-mail: bodrenko@mail.ru \qquad \qquad http://www.bodrenko.com}
\end{center}

\begin{center}
{\bf Abstract}\\
\end{center}

{\small Author of this article created for the first time the
method for finding solutions of the Minkowski problem for closed
surfaces in Riemannian space. }

\section*{Introduction}

\bigskip

Author resolves the Minkowski problem as the problem of
construction the continuous G-deformations preserving the product
of principal curvatures for every point of surface in Riemannian
space. G-deformation transfers every normal vector of surface in
parallel along the path of the translation for each point of
surface.

The Minkowski problem (MP) asks the fundamental question of
differential geometry. H. Minkowski, in 1903, published first
article on this subject.  Generalizations of the MP in Euclidean
and pseudo-Euclidean spaces were made basically by developing the
methods created by H.Minkowski, A.D. Aleksandrov, A.V. Pogorelov,
W.J. Firey. But nobody has studied this problem in Riemannian
space.

The MP in Riemannian space is much more complicate than the MP in
Euclidean space because well known methods that are useful for
 Euclidean and pseudo-Euclidean spaces are not valid for
 Riemannian space.

Author uses term of $G-$deformation because for surface in
Riemannian space there is no such term as spherical image.

\section*{\S 1. Basic definitions. Statement of the main result.}

Let $R^{3}$ be  the three-dimensional Riemannian space with metric
tensor $\tilde{a}_{\alpha\beta},$
  $F$ be the two-dimensional simply connected oriented
closed surface in $R^{3}.$

Let $F\in C^{m,\nu}, \nu \in (0;1) , m\ge 4.$  Let $F$ has all
strictly positive principal curvatures $k_{1}$ and $k_{2}$. Let F
be oriented so that mean curvature $H$ is strictly positive.
Denote $K=k_{1}k_{2}.$

Let $F$ be glued from the two-dimensional simply connected
oriented surfaces $F^{+}$ and $F^{-}$ of class $C^{m,\nu}.$ Let
$F^{+}$ be attached to $F^{-}$ along the common boundary $\partial
F$ of class $C^{m+1,\nu}.$

Let $F^{+}$ and $F^{-}$ be given by immersions of the domain
$D\subset E^{2}$ into $R^{3}$ by the equation:
$y^{\sigma}=f^{\sigma\pm}(x), x\in D$, $f^{\pm}:D \rightarrow
R^{3}.$ Denote by $d\sigma(x)=\sqrt{g} dx^{1}\wedge dx^{2}$ the
area element of the surface $F$. We identify the points of
immersion of surface $F$ with the corresponding coordinate sets in
$R^{3}$. Without loss of generality we assume that $D$ is unit
disk. Let $x^{1}, x^{2}$ be the Cartesian coordinates.

Symbol $_{,i}$ denotes covariant derivative in metric of surface
$F.$ Symbol $\partial_{i}$ denotes partial derivative by variable
$x^{i}.$ We will assume $\dot{f}\equiv \frac{d f}{d t}.$ We define
$\Delta(f)\equiv f(t)-f(0).$ Let $g_{ij}$ and $b_{ij}$ be the
coefficients of the first and the second fundamental form
respectively.

We consider continuous deformation of the surface $F$: $\{F_{t}\}$
defined by the equations
$$y^{\sigma}_{t}=y^{\sigma} + z^{\sigma}(t), z^{\sigma}(0)\equiv 0,
 t\in[0;t_{0}], t_{0}>0. \eqno(1.1)$$

\begin{definition}
\label{definition 1}. Deformation $\{F_{t}\}$ is called the
continuous deformation preserving the product of principal
curvatures  ( or $M-$deformation [21]) if the following condition
holds: $\Delta(K)=0$ and $z^{\sigma}(t)$ is continuous by $t.$

\end{definition}

The deformation $\{F_{t}\}$ generates the following set of paths
in $R^{3}$
$$u^{\alpha_{0}}(\tau)=(y^{\alpha_{0}}+z^{\alpha_{0}}(\tau)), \eqno(1.2)$$
where $z^{\alpha_{0}}(0)\equiv 0, \tau \in [0;t], t\in[0;t_{0}],
t_{0}>0.$

\begin{definition}
\label{definition 2}. The deformation $\{F_{t}\}$ is called the
$G-$deformation if every normal vector of surface transfers in
parallel along the path of the translation for each point of
surface.
\end{definition}

\begin{theorem}
\label{theorem1}. Let $F\in C^{m,\nu}, \nu \in (0;1) , m\ge 4,$ be
closed surface. Let $F$ be glued from the two-dimensional simply
connected oriented surfaces $F^{+}$ and $F^{-}$ of class
$C^{m,\nu}.$ Let $F^{+}$ be attached to $F^{-}$ along the common
boundary $\partial F$ of class $C^{m+1,\nu}.$ Let
$\tilde{a}_{\alpha\beta} \in C^{m,\nu},$ $\exists M_{0}=const>0$
such that $\|\tilde{a}_{\alpha\beta}\|_{m,\nu}<M_{0},$ $\|\partial
\tilde{a}_{\alpha\beta}\|_{m,\nu}<M_{0},$ $\|\partial^{2}
\tilde{a}_{\alpha\beta}\|_{m,\nu}<M_{0}.$

1) Then there exists $t_{0}>0$ such that
 for all $t\in [0, t_{0})$
there exists three-parametric $MG-$deformation of class
$C^{m-2,\nu}$ continuous by $t.$

2) If, at the point $T_{0}\in F^{+},$ the following additional
condition holds: $\forall t : z^{\sigma}(t)\equiv 0.$ Then there
exists $t_{0}>0$ such that
 for all $t\in [0, t_{0})$
there exists  only zero $MG-$deformation of class $C^{m-2,\nu}$
continuous by $t.$

\end{theorem}

\section*{\S 2. Deduction of $G-$deformation formulas for surfaces in
Riemannian space.}

We denote:
$$z^{\sigma}(t)=a^{j}(t)y^{\sigma},_{j} + c(t)n^{\sigma},  \eqno(2.1)$$
where $a^{j}(0)\equiv 0, c(0)\equiv 0,$ $n^{\sigma}$ is unit
normal vector of surface at the point $(y^{\sigma}).$ Therefore
the deformation of surface $F$ is defined by functions $a^{j}$ and
$c.$

 We use designations from [30]. Notice that the
equations of $G-$deformation
obtained in \S2 of [30] are :\\
$$ a^{l}b_{li}+(1+N_{0}(t,0))c,_{i}+
\partial_{i}a^{j} N_{j}(t,0) +Q_{i}(t,0)=0,i=1,2. \eqno(2.1) $$

The estimations of norms for functions $N_{0}(t,0)), N_{j}(t,0),
Q_{i}(t,0) $ are derived in \S 3 of article [30] where the
explicit forms of $N_{0}(t,0)), N_{j}(t,0), Q_{i}(t,0)$ are found.

We introduce conjugate isothermal coordinate system where
$b_{ii}=V,i=1,2, b_{12}=b_{21}=0.$ Then we have the equation
system from (2.1):
$$ c,_{1}(1+N_{0})+Va^{1}+N_{k}\partial_{1}a^{k}+Q_{1}=0$$
$$ c,_{2}(1+N_{0})+Va^{2}+N_{k}\partial_{2}a^{k}+Q_{2}=0 \eqno(2.2)$$

We find solution of the equation system (2.2) by
 finding function $\dot{c}$ on functions $\dot{a}^{i}.$

We denote
$$
\Psi_{1}= -(c_{,1}\partial_{2}N_{0}-c_{,2}\partial_{1}N_{0}+
\partial_{1}a^{k}\partial_{2}N_{k}-
\partial_{2}a^{k}\partial_{1}N_{k}+
\partial_{2}Q_{1}-\partial_{1}Q_{2})/V. \eqno(2.3)
$$
Then, from (2.2) and (2.3), we have the following equation (see
[30]):
$$ \partial_{2}\dot{a}^{1}-\partial_{1}\dot{a}^{2}+
p_{k}\dot{a}^{k}=\dot{\Psi}_{1}, \eqno(2.4) $$ where
$p_{1}=\partial_{2}(\ln V), p_{2}=-\partial_{1}(\ln V).$ Note that
$p_{k}$ do not depend on $t.$

At first, we will solve the equation system (2.2) for surface
$F^{+}$ for the case 2) of theorem 1, assuming that functions
$a^{1}$ and $a^{2}$ are given. The solution of the equation system
(2.2) was presented in [30].

From [30] we have that every pair of functions $\dot{a}^{i}\in
C^{m-2,\nu}$ corresponds to the unique function $\dot{c}\in
C^{m-2,\nu}$ and therefore to the unique function $c \in
C^{m-2,\nu}.$

We find function $\dot{c}$ from the equation:
$$\dot{c}=L_{a}(\dot{c})+\gamma_{t}, \eqno(2.5) $$
by the method of successive approximations. Notice that $L_{a}$
and $\gamma_{t}$ were presented in [30] and have explicit forms.

In [30] was proved the following lemma for the case 2) of theorem
1:

 {\bf Lemma 5.1. ([30])} {\it Let the following conditions hold:

1) metric tensor in $R^{3}$ satisfies the conditions: $\exists
M_{0}=const>0$ such that
$\|\tilde{a}_{\alpha\beta}\|_{m,\nu}<M_{0},$ $\|\partial
\tilde{a}_{\alpha\beta}\|_{m,\nu}<M_{0},$ $\|\partial^{2}
\tilde{a}_{\alpha\beta}\|_{m,\nu}<M_{0}.$

2) $\exists t_{0}>0$ such that $a^{k}(t), \partial_{i}a^{k}(t),
\dot{a}^{k}(t), \partial_{i}\dot{a}^{k}(t)$ are continuous by $t,
\forall t\in[0,t_{0}],$ $ a^{k}(0)\equiv 0,
\partial_{i}a^{k}(0)\equiv 0.$

3) $\exists t_{0}>0$ such that $a^{i}(t)\in C^{m-2,\nu} ,
\partial_{k} a^{i}(t)\in C^{m-3,\nu},$ $\forall t\in[0,t_{0}].$

Then $\exists t_{*}>0$ such that the equation
 $\dot{c}=L_{a}(\dot{c})+\gamma_{t}$ $\forall t\in[0,t_{*}].$
 has unique solution of class $C^{m-2,\nu}$ continuous by $t.$
}

Proof of lemma 5.1. was presented in [30].

{\bf Corollary.([30])}
{ \it Let the conditions of lemma 5.1. hold.\\
Then the function $\dot{c}$ takes the form:
$$ \dot{c}(x^{1},x^{2},t)=
\int\limits_{(x^{1}_{(0)},x^{2}_{(0)})}^{(x^{1},x^{2})}
\Biggl(-V\dot{a}^{1}\Biggr) d\tilde{x}^{1}+
\Biggl(-V\dot{a}^{2}\Biggr)d\tilde{x}^{2}+P(\dot{a}^{1},\dot{a}^{2}),
$$ and for $P$ the following inequality holds:
$$\|P(\dot{a}^{1}_{(1)},\dot{a}^{2}_{(1)})-
P(\dot{a}^{1}_{(2)},\dot{a}^{2}_{(2)})\|_{m-2,\nu}\leq
K_{8}(t)(\|\dot{a}^{1}_{(1)}-\dot{a}^{1}_{(2)}\|_{m-2,\nu}+
\|\dot{a}^{2}_{(1)}-\dot{a}^{2}_{(2)}\|_{m-2,\nu}),$$ where for
any $\varepsilon>0$ there exists $t_{0}>0$ such that
 for all $t\in[0,t_{0})$ the following inequality holds:
$K_{8}(t)<\varepsilon.$ }

The proof follows from construction of function $\dot{c}$ (see
[30]).

{\bf Note.} {\it For surface $F^{-}$ we use similar methods of
finding solutions of the equation system (2.2). But formula (5.6)
from [30] takes the following form
$$ \dot{c}(x^{1},x^{2},t)=$$
$$\int\limits_{(x^{1}_{(0)},x^{2}_{(0)})}^{(x^{1},x^{2})}
\Biggl(- \frac{-V\dot{a}^{1}N_{0}+N_{k}\partial_{1}\dot{a}^{k}+
\dot{N}_{k}\partial_{1}a^{k}+\dot{Q}_{1}}{(1+N_{0})}+ \frac{
\dot{N_{0}}(Va^{1}+N_{k}\partial_{1}a^{k}+Q_{1})}{(1+N_{0})^{2}}
\Biggr) d\tilde{x}^{1}+$$
$$
\Biggl(-\frac{-V\dot{a}^{2}N_{0}+N_{k}\partial_{2}\dot{a}^{k}+
\dot{N}_{k}\partial_{2}a^{k}+\dot{Q}_{2}}{(1+N_{0})}+ \frac{
\dot{N_{0}}(Va^{2}+N_{k}\partial_{2}a^{k}+Q_{2})}{(1+N_{0})^{2}}
\Biggr)d\tilde{x}^{2}+$$
$$\int\limits_{(x^{1}_{(0)},x^{2}_{(0)})}^{(x^{1},x^{2})}
\Biggl(-V\dot{a}^{1}\Biggr) d\tilde{x}^{1}+
\Biggl(-V\dot{a}^{2}\Biggr)d\tilde{x}^{2}+B^{-}_{1}, \eqno(2.6)$$
where $B^{-}_{1}$ is arbitrary real parameter. $B^{-}_{1}$ is
found from the condition $\dot{c}^{+}=\dot{c}^{-}$ on $\partial
F.$}

\section*{\S 3. The equations of $MG-$deformations.}

Deduction the formulas of deformations preserving the product of
principal curvatures was presented in [30]. The equation of
deformations preserving the product of principal curvatures takes
the following form:
$$\partial_{1}\dot{a}^{1}+\partial_{2}\dot{a}^{2}+q^{(b)}_{k}\dot{a}^{k}=
\dot{\Psi}_{2}^{(b)},  \eqno(3.1)$$ where
$\dot{\Psi}_{2}^{(b)}=q^{(b)}_{0}\dot{c}- P_{0}(\dot{a}^{1},
\dot{a}^{2},\partial_{i}\dot{a}^{j}).$ $P_{0}$ has explicit form.
Notice that $q^{(b)}_{k}\in C^{m-3,\nu},$ $q^{(b)}_{0} \in
C^{m-3,\nu}$ and
 do not depend on $t.$

We will use the following lemma form [30].

 {\bf Lemma 6.2.1. ([30])} {\it Let the following conditions hold:

1) metric tensor in $R^{3}$ satisfies the conditions: $\exists
M_{0}=const>0$ such that
$\|\tilde{a}_{\alpha\beta}\|_{m,\nu}<M_{0},$ $\|\partial
\tilde{a}_{\alpha\beta}\|_{m,\nu}<M_{0},$ $\|\partial^{2}
\tilde{a}_{\alpha\beta}\|_{m,\nu}<M_{0}.$

2) $\exists t_{0}>0$ such that $a^{k}(t), \partial_{i}a^{k}(t),
\dot{a}^{k}(t), \partial_{i}\dot{a}^{k}(t)$ are continuous by $t,
\forall t\in[0,t_{0}],$ $ a^{k}(0)\equiv 0,
\partial_{i}a^{k}(0)\equiv 0.$

3) $\exists t_{0}>0$ such that $a^{i}(t)\in C^{m-2,\nu} ,
\partial_{k} a^{i}(t)\in C^{m-3,\nu},$ $\forall t\in[0,t_{0}].$

Then $\exists t_{*}>0$ such that for all $t\in[0,t_{*})$ $P_{0}\in
C^{m-3,\nu}$ and the following inequality holds:
$$\|P_{0}(\dot{a}^{1}_{(1)},\dot{a}^{2}_{(1)})-
P_{0}(\dot{a}^{1}_{(2)},\dot{a}^{2}_{(2)})\|_{m-2,\nu}\leq
K_{9}(t)(\|\dot{a}^{1}_{(1)}-\dot{a}^{1}_{(2)}\|_{m-1,\nu}+
\|\dot{a}^{2}_{(1)}-\dot{a}^{2}_{(2)}\|_{m-1,\nu}),$$ where for
any $\varepsilon>0$ there exists $t_{0}>0$ such that
 for all $t\in[0,t_{0})$ the following inequality holds:
$K_{9}(t)<\varepsilon.$ }

The proof was presented in [30].

We observe all designations for functions and terms presented in
[30]. We use estimations of norms for obtained functions from \S2,
\S3, \S7 of [30]. Formulas of functions $\dot{W}_{1},\dot{W}_{2},
\dot{\Psi}_{2},$ the estimations for norms of  these functions are
presented in \S8 of [30]. From article [30], we also use the
following lemmas: 2.1, 2.2, 2.3, 3.1, 3.2, 5.1, 6.2.1, 7.1, 7.2,
7.3, 8.3.1, 8.3.2, and theorems: 1, 9.1, 9.2.

We have the following equation system of elliptic type from [30]:
$$ \partial_{2}\dot{a}^{1}-\partial_{1}\dot{a}^{2}+p_{k}\dot{a}^{k}=
\dot{\Psi_{1}}, $$
$$ \partial_{1}\dot{a}^{1}+\partial_{2}\dot{a}^{2}+q^{(b)}_{k}\dot{a}^{k}=
\dot{\Psi}^{(b)}_{2}, \eqno(3.2) $$ where we use (3.1).
$\dot{\Psi}^{(b)}_{2}=q^{(b)}_{0}\dot{c}-P_{0}.$
 Note that $q^{(b)}_{k}$ do not depend on $t.$

Then we finally have the form of desired equation from [30]:
$$\partial_{\bar{z}}\dot{w}+A\dot{w}+B\bar{\dot{w}}+E(\dot{w})=\dot{\Psi}. \eqno(3.3) $$

Let, along the $\partial F$, be given vector field tangent to $F.$
We denote it by the following formula:
$$v^{\alpha}=l^{i}y^{\alpha}_{,i}.  \eqno(3.4)$$

We consider the boundary-value condition:
$$\tilde{a}_{\alpha\beta}z^{\alpha}v^{\beta}=\tilde{\gamma}(s,t) ,
s\in \partial D. \eqno(3.5)$$

Define:
$\tilde{\lambda}_{k}=\tilde{a}_{\alpha\beta}y^{\alpha}_{,k}v^{\beta},
k=1,2.$ \\ Then boundary condition takes the form:
$Re\{(\dot{a^{1}}+i \dot{a^{2}}) (\tilde{\lambda}_{1}-i
\tilde{\lambda}_{2})\}=\dot{\tilde{\gamma}}$ on $\partial F.$

Denote: $\lambda_{k}=\frac{ \tilde{\lambda}_{k}}
{(\tilde{\lambda}_{1})^{2}+(\tilde{\lambda}_{2})^{2}}, k=1,2.$
$\dot{\varphi}=\frac{ \dot{\tilde{\gamma}} }
{(\tilde{\lambda}_{1})^{2}+(\tilde{\lambda}_{2})^{2}}.$

Then boundary-value condition takes the form:
$Re\{\overline{\lambda} \dot{w}\}=\dot{\varphi}$ on $\partial F,$
where $|\lambda|=1.$

We analyze the decidability of the following equation (A):\\
$$\partial_{\bar{z}}\dot{w}+A\dot{w}+B\bar{\dot{w}}+E(\dot{w})=\dot{\Psi},
\qquad Re\{\overline{\lambda}\dot{w}\}=\dot{\varphi} \quad on
\quad \partial D, \eqno(3.6) $$
$\lambda=\lambda_{1}+i\lambda_{2},$ $|\lambda|\equiv 1,$ $\lambda,
\dot{\varphi}\in C^{m-2,\nu}(\partial D).$

We will use the fact that $\dot{\Psi}=\dot{\Psi}(\dot{w},z,t),
E(\dot{w})=E(\dot{w},z,t),$
$\dot{w}=\dot{w}(t),$ \\
$\dot{\varphi}=\dot{\varphi}(s,t), s\in \partial D,$
$\lambda=\lambda(s), s\in \partial D.$

Let $n$ be index of obtained boundary-value problem
$$n=\frac{1}{2\pi}\Delta_{\partial D} \arg \lambda(s). \eqno(3.7) $$

Since $v^{\alpha}$ is not tangential
 to $\partial F$ vector field then the index of boundary-value problem
 $n=1.$

\section*{\S 4. Proof of theorem 1.}

{\bf Proof of theorem 1.}

At first, we will prove the case 2) of theorem 1.

 Notice that we consider closed surface. We assume that surface
consists of two surfaces $F^{+}$ and $F^{-}$ with the same
boundary $\partial F.$ Let $F^{+}$ and $F^{-}$ be glued. According
to [17, 18], we introduce conjugate isothermal coordinate system
on $F.$ Therefore $F$ is mapped gomeomorphically onto the domains
$G^{+}\subset E^{2}$ and $G^{-}\subset E^{2}.$
 $G^{+}$ and $G^{-}$ have the same boundary $\Gamma$ and
 $G^{-}$ is infinite domain. Notice that $\Gamma$ is simple
 closed smooth curve of class $C^{m+1,\nu}$.
Note that $\partial F$ is mapped gomeomorphically onto the curve
$\Gamma$
 twice.
According to [17, 18], without loss of generality we assume that
for each
 point $M\in \partial F$ there exists unique point $s_{M}\in \Gamma.$

It is possible to introduce uniform conjugate isothermal
coordinate
 system on $F$ that maps $F$ gomeomorphically onto the whole $E^{2}$
 if we do not fix domains $G^{+}$ and $G^{-}$ ( see [17, 18] ).

Along the $\partial F,$ we consider the condition of contiguity of
surfaces $F^{+}$ and $F^{-}.$ Let
$s^{\alpha}=s^{k}y^{\alpha}_{,k}$ be unit tangential to the
$\partial F$ vector field (see [28]). Let $v^{\alpha}$ be unit
vector field which, on $\partial F$, satisfies the following
conditions:
$$v^{\alpha}=l^{i}y^{\alpha}_{,i}, \quad
\tilde{a}_{\alpha\beta}v^{\alpha}s^{\beta}=0, \quad
\tilde{a}_{\alpha\beta}v^{\alpha}v^{\beta}=1.$$ Denote:
$z^{+\alpha}_{s}=\tilde{a}_{\alpha\beta}z^{+\alpha}s^{\beta},$
$z^{-\alpha}_{s}=\tilde{a}_{\alpha\beta}z^{-\alpha}s^{\beta},$
$z^{+\alpha}_{v}=\tilde{a}_{\alpha\beta}z^{+\alpha}v^{\beta},$
$z^{-\alpha}_{v}=\tilde{a}_{\alpha\beta}z^{-\alpha}v^{\beta}.$

Therefore, according to [17,18], if the surfaces $F^{+}$ and
$F^{-}$ compose the smooth surface $F$ then the condition of
contiguity of surfaces $F^{+}$ and $F^{-}$ takes the form:
$$z^{+\alpha}_{s}=z^{-\alpha}_{s}, \quad z^{+\alpha}_{v}=z^{-\alpha}_{v},
\quad c^{+}=c^{-} \quad on \quad \partial F. $$

We have the following conditions:\\
$$\tilde{a}_{\alpha\beta}v^{\alpha}s^{\beta}=0,$$
$$\tilde{a}_{\alpha\beta}v^{\alpha}v^{\beta}=1,$$
$$\tilde{a}_{\alpha\beta}s^{\alpha}s^{\beta}=1.$$

Therefore we get the following equation system:\\
$$g_{ik}l^{i}s^{k}=0.$$
$$g_{ik}l^{i}l^{k}=1.$$
$$g_{ik}s^{i}s^{k}=1.$$

Consider the formula:\\
$$z^{\alpha}=a^{i} y^{\alpha}_{,i}+c n^{\alpha}. $$

Then
$$z_{v}=g_{ik}a^{i}l^{k},$$
$$z_{s}=g_{ik}a^{i}s^{k}.$$

Therefore we have the following boundary-value condition on $\partial F$:\\
$$a^{+i}=a^{-i}.$$

We obtain on $\partial F$:\\
$$\dot{a}^{+i}=\dot{a}^{-i}.$$

In the article [30]
author analyzed the decidability of the following boundary-value problem (A):\\
$$\partial_{\bar{z}}\dot{w}+A\dot{w}+B\bar{\dot{w}}+E(\dot{w})=\dot{\Psi},
\qquad Re\{\overline{\lambda}\dot{w}\}=\dot{\varphi} \quad on
\quad \partial D, \eqno(4.1) $$ $|\lambda|\equiv 1,$
$\lambda=\lambda_{1}+i\lambda_{2},$ $\lambda, \dot{\varphi}\in
C^{m-2,\nu}(\partial D).$

Therefore we have the following boundary-value problem $(A')$:\\
$$\partial_{\bar{z}}\dot{w}+A\dot{w}+B\bar{\dot{w}}+E(\dot{w})=\dot{\Psi},
\quad (in \quad G^{+} \quad and \quad G^{-}), $$

$$  \dot{a}^{+i}=\dot{a}^{-i}, \quad on \quad \Gamma. \eqno(4.2) $$

The boundary-value problem $(A')$ takes the form:\\
$$\partial_{\bar{z}}\dot{w}+A\dot{w}+B\bar{\dot{w}}+E(\dot{w})=\dot{\Psi},
\quad (in \quad G^{+} \quad and \quad G^{-}), $$

$$  \dot{w}^{+}=\dot{w}^{-}, \quad
on \quad \Gamma.  \eqno(4.3) $$

The above problem is generalized Carleman problem. Since $F\in
C^{m,\nu}, m \ge 4$ then the condition of contiguity of surfaces
$F^{+}$ and $F^{-}$ is of at least 4-th order.

Consider the following problem
$$\partial_{\bar{z}}\dot{w}+A\dot{w}+B\bar{\dot{w}}=0,
\quad (in \quad G^{+} \quad and \quad G^{-}), $$
$$  \dot{w}^{+}=\dot{w}^{-}, \quad
on \quad \Gamma.  \eqno(4.4) $$

Using methods from [17,18] we reduce boundary-value problem (4.4)
to the following problem for analytic functions: find analytic
function $\Phi^{+}(z)$ in $G^{+}$ and analytic function
$\Phi^{-}(z)$ in $G^{-}$ satisfying the boundary-value condition
$\Phi^{+}(s)$=$U(s)\Phi^{-}(s)$ on $\Gamma,$ where $U(s)\neq 0$ on
$\Gamma,$ index $\kappa=$ ind $U(s)=0,$  $\Phi^{-}$ has finite
order on infinity $\Phi^{-}(\infty)=0$ and $\Phi^{+}(z_{0})=0$ at
the point $z_{0}=x^{1}_{0}+ix^{2}_{0}\in G^{+}.$ Therefore from
[17, 18, 31] we obtain that the boundary-value problem (4.4) has
only zero solution of class $C^{m-2,\nu}.$ $(\Phi^{+}\equiv 0,$
$\Phi^{-}\equiv 0).$

For every admissible analytic in $G^{+}$ function $\Phi^{+}$ the
equation
$$\partial_{\bar{z}}\dot{w}+A\dot{w}+B\bar{\dot{w}}+E(\dot{w})=\dot{\Psi},
\quad in \quad G^{+} $$ has unique solution and is solved by
method of successive approximations as (see [17, 18]):
$$
\dot{w}_{0}(z)=\Phi^{+}(z)+T(-E(\dot{w})+\dot{\Psi}), $$
$$
\dot{w}_{k+1}(z)=\frac{1}{\pi}\int\int\limits_{G^{+}}\frac{(A\dot{w}_{k}+B\overline{\dot{w}}_{k})}{\zeta-z}d\xi
d\eta+\Phi^{+}(z)+T(-E(\dot{w})+\dot{\Psi}), \quad k=0,1,..., $$
where $\zeta=\xi+i \eta,$
$$T(f)=-\frac{1}{\pi}\int\int\limits_{G^{+}}\frac{f(\zeta)}{\zeta-z}d\xi d\eta.$$
Then we have the following inequality: $\|\dot{w}\|_{m-2,\nu}\leq
K_{20}(t)\|\dot{w}\|_{m-2,\nu}$ for all sufficiently small $t\geq
0,$ where $K_{20}(t)<1.$

Using methods from [17, 18, 31, 20] for finding the solutions of
problem (4.3), theory of completely continuous operators and
theory of Volterra operator equation, the theorem 1 from [30] we
obtain that the problem $(A')$ has only zero solution of class
$C^{m-2,\nu}$ continuous by $t.$

The case 2) of theorem 1 is proved.

Proof of case 1) differs from proof of case 2) by the fact that
$\dot{c}^{+}$ has one additional real parameter and $\Phi^{+}$ is
two-parametric function (with two real parameters).

\section*{\bf References.}

\begin{enumerate}

 \item A.I. Bodrenko.
On continuous almost ARG-deformations of hypersurfaces in
Euclidean space [in Russian]. Dep. in VINITI 27.10.1992.,
N3084-T92, UDK 513.81, 14 pp.

 \item A.I. Bodrenko. Some properties continuous ARG-deformations
 [in Russian]. Theses of international science conference
"Lobachevskii and modern geometry", Kazan, Kazan university
publishing house, 1992 ., pp.15-16.

 \item A.I. Bodrenko. On continuous ARG-deformations [in Russian].
Theses of reports on republican science and methodical conference,
dedicated to the 200-th anniversary of N.I.Lobachevskii, Odessa,
Odessa university publishing house, 1992 ., Part 1, pp.56-57.

 \item A.I. Bodrenko. On extension of infinitesimal almost
 ARG-deformations closed \\
 hypersurfaces into analytic
 deformations in Euclidean spaces [in Russian].
 Dep. in VINITI 15.03.1993., N2419-T93 UDK 513.81, 30 pp.

 \item A.I. Bodrenko. On extension of infinitesimal almost
 ARG-deformations of hypersurface with boundary into analytic deformations
[in Russian]. Collection works of young scholars of VolSU,
Volgograd, Volgograd State University publishing house, 1993,\\
pp.79-80.

 \item A.I. Bodrenko. Some properties of continuous almost
AR-deformations of hypersurfaces with prescribed change of
Grassmannian image  [in Russian].
 Collection of science works of young scholars, Taganrog,
 Taganrog State Pedagogical Institute publishing house , 1994, pp. 113-120.

 \item A.I. Bodrenko. On continuous almost AR-deformations
with prescribed change of Grassmannian image [in Russian].
All-Russian school-colloquium on stochastic \\ methods of geometry
and analysis. Abrau-Durso. Publisher Moscow: "TVP". Theses of
reports, 1994, pp. 15-16.

 \item A.I. Bodrenko. Extension of infinitesimal almost ARG-deformations
of \\ hypersurfaces into analytic deformations [in Russian].
All-Russian school-colloquium on stochastic methods. Yoshkar-Ola.
Publisher Moscow: "TVP". Theses of reports, 1995, pp. 24-25.

 \item A.I. Bodrenko. Areal-recurrent deformations of hypersurfaces
 preserving Grassmannian image [in Russian].
Dissertation of candidate of physical-mathematical sciences. \\
Novosibirsk,1995, pp. 85.

 \item A.I. Bodrenko. Areal-recurrent deformations of hypersurfaces
preserving Grassmannian image [in Russian]. Author's summary of
dissertation of candidate of physical- \\ mathematical sciences.
Novosibirsk, 1995, pp. 1-14.

 \item A.I. Bodrenko. Some properties of ARG-deformations [in Russian].
Izvestiay Vuzov. Ser. Mathematics, 1996, N2, pp.16-19.

 \item A.I. Bodrenko. Continuous almost ARG-deformations of surfaces
 with boundary [in Russian].
Modern geometry and theory of physical fields. \\ International
geometry seminar of N.I.Lobachevskii Theses of reports, Kazan,
\\
Publisher Kazan university, 1997, pp.20-21.

 \item A.I. Bodrenko. Continuous almost AR-deformations of surfaces
with prescribed change of Grassmannian image [in Russian]. Red.
"Sib. mat. zhurnal.", Sib. otd. RAN , Novosibirsk, Dep. in VINITI
13.04.1998., N1075-T98 UDK 513.81, 13 pp.

 \item A.I. Bodrenko. Almost ARG-deformations of the second order of
 surfaces in \\ Riemannian space [in Russian].
 Surveys in Applied and Industrial Mathematics. \\ 1998,
 Vol. 5, Issue 2, p.202. Publisher Moscow: "TVP".

 \item A.I. Bodrenko. Almost $AR$-deformations of a surfaces with prescribed
change of \\ Grassmannian image with exterior connections [in
 Russian]. Red. zhurn. "Izvestya vuzov. Mathematics.", Kazan, Dep.
in VINITI  03.08.1998, N2471 - B 98. P. 1-9.

 \item A.I. Bodrenko. Properties of generalized G-deformations with
areal condition of normal type in Riemannian space [in Russian].
 Surveys in Applied and Industrial Mathematics.
Vol. 7. Issue 2. (VII All-Russian school-colloquium on stochastic
methods. Theses of reports.) P. 478. Moscow: TVP, 2000.

 \item I.N. Vekua. Generalized Analytic Functions. Pergamon.
 New York. 1962.

 \item I.N. Vekua. Generalized Analytic Functions [in Russian].
 Moscow. Nauka. 1988.

 \item I.N. Vekua. Some questions of the theory of differential equations
 and  applications in mechanics [in Russian].
 Moscow:"Nauka". 1991 . pp. 256.

 \item A.V. Zabeglov. On decidability of one nonlinear
  boundary-value problem for AG-deformations
  of surfaces with boundary [in Russian].
  Collection of science works. \\ Transformations of surfaces,
  Riemannian spaces determined by given recurrent \\ relations. Part 1.
Taganrog. Taganrog State Pedagogical Institute publishing house.
1999. pp. 27-37.

 \item V.T. Fomenko. On solution of the generalized Minkowski problem
 for surface with boundary [in Russian].
Collection of science works. Transformations of surfaces,\\
  Riemannian spaces determined by given recurrent relations. Part 1.
Taganrog. \\ Taganrog State Pedagogical Institute publishing
house. 1999. pp. 56-65.

 \item V.T. Fomenko. On uniqueness of solution of the generalized
 Christoffel problem for surfaces with boundary [in Russian].
Collection of science works. Transformations of surfaces,
  Riemannian spaces determined by given recurrent relations. Part 1.
Taganrog. Taganrog State Pedagogical Institute publishing house.
1999. pp. 66-72.

 \item V.T. Fomenko. On rigidity of surfaces with boundary in
 Riemannian space [in Russian].
 Doklady Akad. Nauk SSSR. 1969 . Vol. 187, N 2, pp. 280-283.

 \item V.T. Fomenko. ARG-deformations of hypersurfaces in Riemannian
space [in Russian].\\  //Dep. in VINITI 16.11.1990 N5805-B90

  \item S.B. Klimentov. On one method of construction the solutions
  of boundary-value \\ problems in the bending theory of surfaces of positive
  curvature [in Russian]. \\ Ukrainian geometry sbornik. pp. 56-82.

 \item M.A. Krasnoselskii. Topological methods in the theory of
 nonlinear problems

    [in Russian]. Moscow, 1965.

 \item L.P. Eisenhart. Riemannian geometry [in Russian]. \\ Izd. in. lit.,
 Moscow 1948.
 (Eisenhart Luther Pfahler. Riemannian geometry. 1926.)

 \item J.A. Schouten, D.J. Struik. Introduction into new methods
 of differential geometry [in Russian]. Volume 2. Moscow. GIIL. 1948 .
 (von J.A. Schouten und D.J. Struik.
 Einf$\ddot{u}$hrung in die neueren methoden der \\
 differentialgeometrie.
 Zweite vollst$\ddot{a}$ndig umgearbeitete
 Auflage. Zweiter band. 1938. )

 \item I. Kh. Sabitov. //VINITI. Results of science and technics.
 Modern problems of \\mathematics [in Russian].
 Fundamental directions. Vol.48, pp.196-271.

\item Andrei I. Bodrenko. The solution of the Minkowski problem
for open surfaces in Riemannian space. Preprint 2007.
arXiv:0708.3929v1 [math.DG]

\item F.D. Gakhov. Boundary-value problems [in Russian]. GIFML.
Moscow. 1963.

\end{enumerate}

\end{document}